\documentclass{article}

\usepackage{amsmath}    
\usepackage{graphicx} 
\usepackage{verbatim} 
\usepackage{color} 
\usepackage{subfigure}  
\usepackage{hyperref}  
\usepackage{fullpage}
\usepackage{bbm}
\usepackage{dsfont} 
\usepackage{amssymb}
\usepackage{mathrsfs}
\usepackage{enumerate}
\usepackage[super]{nth}
\usepackage{adjustbox}

\usepackage{moreverb}
\usepackage{algorithm}
\usepackage{algorithmicx}
\usepackage{algpseudocode}
\usepackage{pgfplots}
\usepackage{alltt}
\usepackage{listings}
\usepackage{xspace}
\pgfplotsset{compat=1.12}

\newcommand{\sgnode}[1]{{\bf \left<#1\right>}}
\newcommand{\gr}[1]{{\color{gray} #1}}

\newtheorem{theorem}{Theorem}
\newtheorem{prop}{Proposition}
\newtheorem{example}{Example}
\newtheorem{conjecture}{Conjecture}
\newtheorem{lemma}{Lemma}

\newtheorem{definition}{Definition}

\newcommand{\Z}{\mathbb{Z}}
\newcommand{\R}{\mathbb{R}}

\newcommand{\N}{\mathbb{N}}

\renewcommand{\phi}{\varphi}

\newcommand{\Ap}{{\operatorname{Ap}}}

\begin{document}

\title{Counting Numerical Semigroups}
\markright{Counting Semigroups}
\author{Nathan Kaplan}

\maketitle

\begin{abstract}
A numerical semigroup is an additive submonoid of the natural numbers with finite complement.  The size of the complement is called the genus of the semigroup.  How many numerical semigroups have genus equal to $g$?  We outline Zhai's proof of a conjecture of Bras-Amor\'os that this sequence has Fibonacci-like growth.  We now know that this sequence asymptotically grows as fast as the Fibonacci numbers, but it is still not known whether it is nondecreasing.  We discuss this and other open problems.  We highlight the many contributions made by undergraduates to problems in this area.
\end{abstract}

\section{What are we counting?}\label{Intro}

A \emph{numerical semigroup} $S$ is an additive submonoid of $\N_0 = \{0,1,2,\ldots\}$, where $\N_0 \setminus S$ is finite.  We say that $\{n_1,\ldots, n_t\}$ is a \emph{generating set} of $S$ if the elements of $S$ are exactly the linear combinations of $n_1,\ldots, n_t$ with nonnegative integer coefficients.  In this case, we write 
\[
S = \langle n_1, \ldots, n_t \rangle = \left\{a_1 n_1 + \cdots + a_t n_t\ |\ a_1,\ldots, a_t \in \Z_{\ge 0}\right\}.
\] 

A nice way to see that a finite generating set exists is to produce one. The smallest nonzero element of $S$ is called the \emph{multiplicity} of $S$, denoted $m(S)$.  The \emph{Ap\'ery set} of $S$ with respect to $m$ is 
\[
\Ap(m,S) = \{0,k_1 m + 1, k_2 m + 2,\ldots, k_{m-1} m + (m-1)\},
\] 
where $k_1,\ldots, k_{m-1}$ are positive integers defined so that $k_i m + i$ is the smallest positive integer in $S$ that is congruent to $i$ modulo $m$ \cite{Apery}.  It is easy to check that every element of $S$ can be written as a linear combination of elements of $\Ap(m,S)$.  Every numerical semigroup has a unique minimal generating set, which we get by removing any element $k_i m + i$ that is a linear combination of the other Ap\'ery set elements.  The size of this minimal generating set is called the \emph{embedding dimension} of $S$, denoted $e(S)$, and its elements are called \emph{minimal generators}.  For detailed proofs, and an excellent introduction to the subject, see \cite[Chapter 1]{GarciaSanchezRosales}.

Many mathematicians first encounter numerical semigroups through the \emph{linear Diophantine problem of Frobenius} or \emph{Frobenius problem}, which asks for a formula in terms of the minimal generating set for the largest element of the complement $\N \setminus S$.  
\begin{definition}
Let $S$ be a numerical semigroup.  The elements of the complement $\N \setminus S$ are called the \emph{gaps} of $S$.  The largest of these gaps is called the \emph{Frobenius number} of $S$, denoted $F(S)$.  The number of gaps is called the \emph{genus} of $S$, denoted $g(S)$.
\end{definition}

\begin{prop}[Sylvester]
Let $a < b$ be relatively prime positive integers and $S = \langle a,b\rangle$. Then
\begin{enumerate}
\item $F(S) = ab-a-b$,
\item $g(S) = \frac{(a-1)(b-1)}{2}$.
\end{enumerate}
\end{prop}

Selmer observed that the Frobenius number and genus can easily be deduced from $\Ap(m,S)$ \cite{Selmer}.  More specifically, $F(S) = \max \Ap(m,S) - m$ and $g(S) = \sum_{i=1}^{m-1} k_i$.  For $S = \langle a,b\rangle,\ \Ap(a,S) = \{0,b,2b,\ldots, (a-1)b\}$, which implies $F(S) = (a-1)b - a$.  It is not difficult to check that for every $x< F(S)$ exactly one of $\{x,F(S) - x\}$ is contained in $S$, which proves the second statement.  The Frobenius problem for numerical semigroups with three or more generators is an active area of research.  See the book \cite{RamirezAlfonsin} for an excellent overview.

We would like to understand the infinite set of all numerical semigroups, which means that we need a way to order them.  Let $N(g)$ be the number of numerical semigroups $S$ with $g(S) = g$.  Understanding this sequence is the main goal of this article.
\begin{example}
\begin{enumerate}
\item $N(0) = 1$:  The unique numerical semigroup of genus $0$ is $\N_0$.

\item $N(1) = 1$: The only numerical semigroup containing $1$ is $\N_0$, so if $g(S) = 1$ then $\N \setminus S = \{1\}$, which implies $S = \langle 2,3\rangle$.

\item $N(2) = 2$: If $S$ is a numerical semigroup of genus $2$, then $\N \setminus S$ consists of $1$ and exactly one other element.  This second gap must be either $2$ or $3$, because the only numerical semigroups containing both $2$ and $3$ are $\langle 2,3\rangle$ and $\N_0$.  If this second gap is $2$, then $S = \langle 3,4,5\rangle$ and if it is $3$, then $S = \langle 2,5\rangle$.

\end{enumerate}
\end{example}
\noindent Computing these values quickly becomes too complicated to do by hand.  See Figure \ref{Fig1} for some small values of $N(g)$.

\begin{figure}
\centering \scalebox{.84}{
\begin{tabular}{| c | c c c c c c c c c c c c c c c c|}
\hline
$g$ & 0 & 1 & 2 & 3&4&5&6&7&8&9&10 &11 &12 &13 &14 &15 \\
\hline
$N(g)$&1&1&2&4&7&12&23&39&67&118&204&343 &592 &1001 &1693 &2857 \\
\hline
\end{tabular}}
\caption{The number of numerical semigroups of genus $g$ for $g\le 15$.}
\label{Fig1}
\end{figure}

In 2008, Maria Bras-Amor\'os computed $N(g)$ for $g\le 50$ and noticed some striking patterns \cite{BrasAmorosFibonacci}. 
\begin{conjecture}\label{BrasAmorosConjecture}[Bras-Amor\'os]
\begin{enumerate}
\item $N(g) \ge N(g-1) + N(g-2)$, for $g \ge 2$,
\item $\lim_{g \to \infty} \frac{N(g-1) + N(g-2)}{N(g)} = 1$,
\item $\lim_{g\to \infty} \frac{N(g)}{N(g-1)} = \phi$, where $\phi = \frac{1+\sqrt{5}}{2}$ is the golden ratio.
\end{enumerate}
\end{conjecture}
\noindent Note that the third statement implies the second.  

\begin{figure}
\centering
\includegraphics[width=110mm]{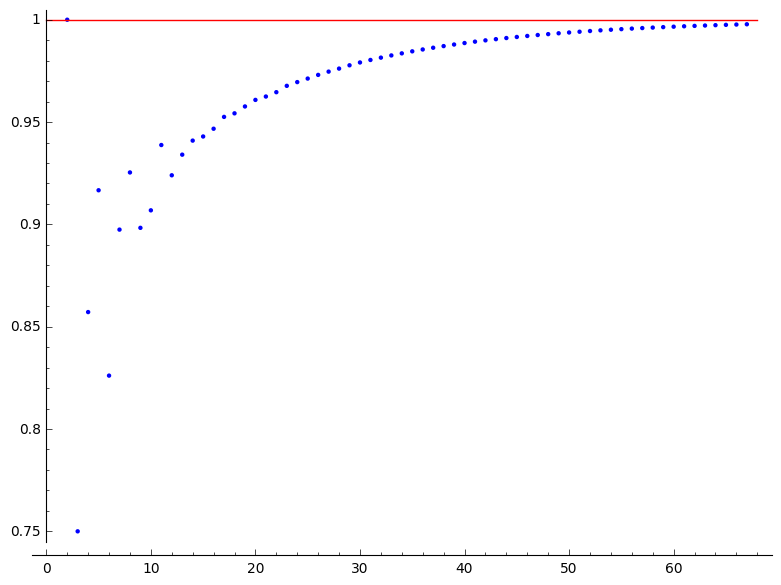}
\caption{Behavior of the quotient $\frac{N(g-1)+N(g-2)}{N(g)}$.}
\label{Fig2}
\end{figure}

\begin{figure}
\centering
\includegraphics[width=110mm]{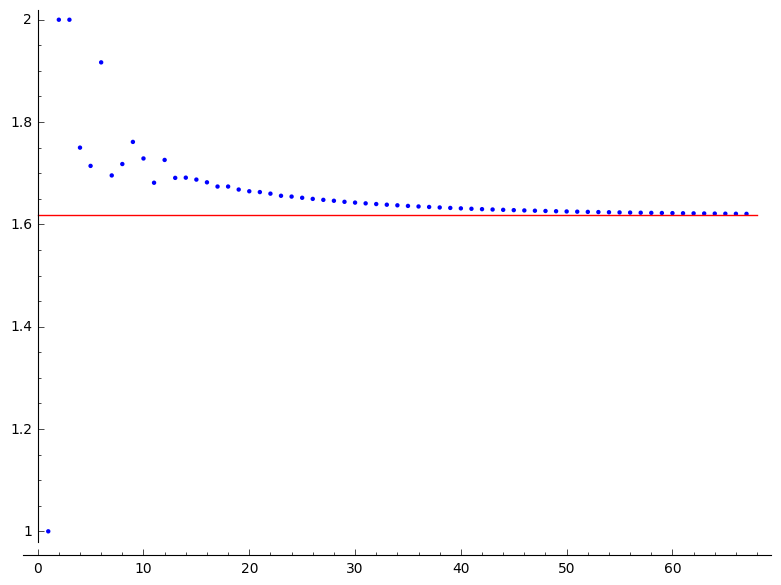}
\caption{Behavior of the quotient $\frac{N(g)}{N(g-1)}$. The horizontal line is $y = \phi$.}
\label{Fig3}
\end{figure}

These computations have been extended by Fromentin and Hivert to $g \le 67$ \cite{FromentinHivert}.  Using similar ideas, Delgado, Garc\'ia-S\'anchez, and Morais have implemented a program to find the set of all numerical semigroups of genus $g$ in the \emph{NumericalSgps} package for the computer algebra system GAP \cite{DelgadoGarciaSanchezMorais, GAP}.  Figures \ref{Fig2} and \ref{Fig3} are updated versions of charts in \cite{BrasAmorosFibonacci} that give computational evidence for Conjecture \ref{BrasAmorosConjecture}.

Alex Zhai proved parts 2 and 3 of Conjecture \ref{BrasAmorosConjecture} while he was an undergraduate \cite{Zhai}.
\begin{theorem}[Zhai]\label{ZhaiTheorem}
Let $N(g)$ be the number of numerical semigroups of genus $g$. Then
\[
\lim_{g\to \infty} \frac{N(g)}{\phi^g} = S,
\]
where $S$ is a constant that is at least $3.78$.
\end{theorem}

The first part of Conjecture \ref{BrasAmorosConjecture} is completely open.  In fact, the much weaker conjecture that $N(g) \ge N(g-1)$ is still unsolved.  Of course, Theorem \ref{ZhaiTheorem} implies that this can fail to hold for only finitely many values of $g$.
\begin{conjecture}\label{IncreasingConjecture}
For all $g \ge 1$ we have $N(g) \ge N(g-1)$.
\end{conjecture}

\subsection{Outline of the paper.}

We first give some potential explanations for the Fibonacci-like growth of $N(g)$.  We then introduce the semigroup tree and describe the strategy for the proof of Theorem \ref{ZhaiTheorem}.  In the final two sections of the paper, we discuss several refined counting questions and other open problems about numerical semigroups.

\section{Why Fibonacci?}\label{Fibonacci}

\subsection{Numerical semigroups with $F(S) < 3m(S)$.}\ \\

We begin by describing some results of Yufei Zhao, then an undergraduate, that give a Fibonacci-like lower bound for $N(g)$ \cite{Zhao}.  Instead of considering all numerical semigroups, he focuses on semigroups where the Frobenius number is bounded in terms of the multiplicity.  The Fibonacci numbers are defined by $F_0 = 0,\ F_1 = 1$, and $F_n = F_{n-1}+F_{n-2}$ for all $n \ge 2$. Recall that $F_n = \frac{1}{\sqrt{5}}\left(\phi^n + (1-\phi)^n\right)$.

\begin{prop}[Zhao]
For any positive integer $g$, the number of numerical semigroups $S$ with genus $g$ satisfying $F(S) < 2 m(S)$ is $F_{g+1}$.
\end{prop}

Zhao first shows that the numerical semigroups with $F(S) < 2 m(S)$ are exactly those consisting of $0, m(S)$, all integers greater than or equal to $2m(S)$, and any subset $A \subseteq [m(S)+1, 2m(S)-1]$.  Such a semigroup satisfies $g(S) = 2(m(S)-1)-|A|$.   Therefore, the total number of numerical semigroups of genus $g$ with $F(S) < 2m(S)$ is
\[
\sum_m \binom{m-1}{2(m-1) - g} = \sum_m \binom{m-1}{g-(m-1)} = F_{g+1}.
\]

Zhao adapts this strategy for semigroups satisfying $2m(S) < F(S) < 3 m(S)$.  
\begin{definition}
For a positive integer $k$, let 
\[
\mathcal{A}_k = \{A \subset [0,k-1]:\ 0 \in A\ \text{ and } k \not\in A + A\},
\]
where $A+A = \{a_1 + a_2:\ a_1, a_2, \in A\}$.

Let $S$ be a numerical semigroup with multiplicity $m$ and Frobenius number $F$, such that $2m < F < 3m$.  We say that $S$ has \emph{type} $(A;k)$, where $k<m$ is a positive integer and $A \in \mathcal{A}_k$, if $F = 2m + k$, and $S \cap [m,m+k] = \{m + a:\ a \in A\}$.
\end{definition}

Zhao shows that the number of numerical semigroups with fixed multiplicity $m(S)$ satisfying $2m(S) < F(S) < 3m(S)$ and type $(A;k)$ can also be expressed as a binomial coefficient, and a similar combinatorial identity gives an expression in terms of Fibonacci numbers.  Let $t(g)$ denote the number of numerical semigroups of genus $g$ satisfying $F(S) < 3m(S)$.
\begin{prop}[Zhao]\label{ZhaoPropB}
For any positive integer $g$, we have
\[
t(g) \ge F_{g+1} + \sum_{k=1}^{\lfloor \frac{g}{3}\rfloor} \sum_{A \in \mathcal{A}_k} F_{g-|(A+A) \cap [0,k]|+ |A|- k - 1}.
\]
\end{prop}
Combining this result with asymptotic estimates for Fibonacci numbers gives a lower bound that is exponential in $\phi$.  Computing the sets $\mathcal{A}_k$ for all $k \le 46$ shows that if $N(g)$ is asymptotic to a constant times $\phi^g$, then that constant is at least $3.78$.

The main step in Zhai's proof of Theorem \ref{ZhaiTheorem} is to prove the following.  
\begin{conjecture}[Zhao]\label{ZhaoConjecture}
We have
\[
\lim_{g\to\infty} \frac{t(g)}{N(g)} = 1.
\]
\end{conjecture}
\noindent We discuss this further in Section \ref{SemigroupTree}.

\subsection{Numerical semigroups with fixed multiplicity.}\ \\

Let $N(m,g)$ be the number of numerical semigroups with multiplicity $m$ and genus $g$.  Figure \ref{Fig4} gives part of a larger table from \cite{Kaplan}.  Kaplan computes these values in an attempt to understand the Fibonacci-like growth of $N(g)$.  A main result of \cite{Kaplan} is that if $m$ and $g$ satisfy certain conditions, then the Fibonacci recurrence holds exactly.

\begin{figure}
\centering
\begin{tabular}{| c || c c c c c c c c c c c  | c |}
\hline

 g$\backslash$m & 1 & 2 & 3 & 4 & 5 & 6 & 7 & 8 & 9 & 10 & 11 &  N(g) \\

\hline
  0& 1 &  &  &  &  & &  &  &  &     &  &  1 \\
 
  1& & 1&  &  &  &  & &  &  &  &    &   1 \\

  2& & 1 & 1&  &  &  &  & &  &  &  &   2 \\
   
  3& & 1 & 2 & 1 &  &  &  &  &  &  &  &  4 \\

  4& & 1 & 2 & 3 & 1 &  &  &  &  &  &  & 7 \\

  5& & 1 & 2 & 4 & 4 & 1 &  &  &  &  &  &   12 \\
      
  6& & 1 & 3 & 6 & 7 & 5 & 1 &  &  &  &  & 23 \\

  7 & & 1 & 3 & 7 & 10 & 11 & 6 & 1 &  &  &  &39 \\

  8 & & 1 & 3 & 9 & 13 & 17 & 16 & 7 & 1 &  & & 67 \\
         
  9 & & 1 & 4 & 11 & 16 & 27 & 28 & 22 & 8 & 1 & &  118 \\

  10 & & 1 & 4 & 13 & 22 & 37 & 44 & 44 & 29 & 9 & 1 &  204 \\
\hline

\end{tabular}
\caption{The number of numerical semigroups of multiplicity $m$ and genus $g$.}
\label{Fig4}
\end{figure}

\begin{theorem}[Kaplan]\label{2g3mTheorem}
Suppose $m$ and $g$ are positive integers satisfying $2g < 3m$.  Then $N(m-1, g-1) + N(m-1,g-2) = N(m,g)$.
\end{theorem}
The proof of this theorem is via an explicit bijection of Ap\'ery sets.  Let $\Ap(m,S) = \{0,k_1 m + 1,\ldots, k_{m-1} m + (m-1)\}$ be the Ap\'ery set of a numerical semigroup with multiplicity $m$.  The \emph{Ap\'ery tuple} or \emph{Kunz coordinate vector} is the tuple of positive integers $(k_1,\ldots, k_{m-1}) \in \Z^{m-1}$, and is a convenient tool for understanding semigroups in terms of their Ap\'ery sets \cite{Kunz}.  When $2g<3m$, Kaplan gives a bijection between Ap\'ery tuples of numerical semigroups with multiplicity $m$ and genus $g$ and Ap\'ery tuples of numerical semigroups with multiplicity $m-1$ and genus $g-1$ or $g-2$.  The bijection takes $(k_1,\ldots, k_{m-1})$ to $(k_1,\ldots, k_{m-2})$.  This result suggests Fibonacci-like growth for $N(g)$ if it is the case that ``most'' numerical semigroups satisfy $2g<3m$, which is not at all clear from the data for small $g$.

\section{The semigroup tree.}\label{SemigroupTree}

Bras-Amor\'os gave some of the earliest upper and lower bounds for $N(g)$ in \cite{BrasAmorosBounds}, showing that 
\[
2F_g \le N(g) \le 1 + 3\cdot 2^{g-3}.
\]
The lower bound is of particular interest because it is asymptotic to a constant times $\phi^g$. This lower bound comes from considering the \emph{semigroup tree} and gives our third explanation for the Fibonacci-like growth of $N(g)$.

The semigroup tree is a rooted tree where the nodes at level $g$ correspond to the numerical semigroups of genus $g$.  Therefore, in order to understand the growth of $N(g)$ we need only understand the number of nodes of each level in the tree.  This is an example of a \emph{Frobenius variety}, an object that nicely organizes families of numerical semigroups closed under certain operations \cite[Chapter 6]{GarciaSanchezRosales}.

The easiest way to specify the semigroup tree is to describe the unique path from any numerical semigroup $S$ back to $\N_0$, the root of the tree.  For a numerical semigroup $S$ of genus $g$ it is easy to check that $S' = S \cup \{F(S)\}$ is a numerical semigroup of genus $g-1$.  Note that $F(S) > F(S')$.  Adjoining $F(S')$ to $S'$ gives a semigroup of genus $g-2$, and we see in this way that starting from $S$ we get a path of $g+1$ semigroups, one of each genus $g' \le g$, ending at $\N_0$.  

Given a numerical semigroup $S$, how many semigroups $S'$ satisfy $S = S' \cup \{F(S')\}$?  The key observation is that for $x > F(S),\ S \setminus \{x\}$ is a numerical semigroup if and only if $x$ is a minimal generator of $S$.  This gives a description of the semigroup tree starting from the root.  
\begin{definition}
Let $S$ be a numerical semigroup of genus $g$.  The \emph{effective generators} of $S$ are the elements of its minimal generating set that are larger than $F(S)$.  The number of effective generators of $S$ is called the \emph{efficacy} of $S$, denoted $h(S)$.

The \emph{children} of a numerical semigroup $S$ of genus $g$ are the numerical semigroups of genus $g+1$ that come from removing an effective generator from $S$.  The number of children of $S$ is $h(S)$.
\end{definition}
%Figure \ref{FigTree}, which is taken from \cite{FromentinHivert}, gives the first few levels of the semigroup tree.

 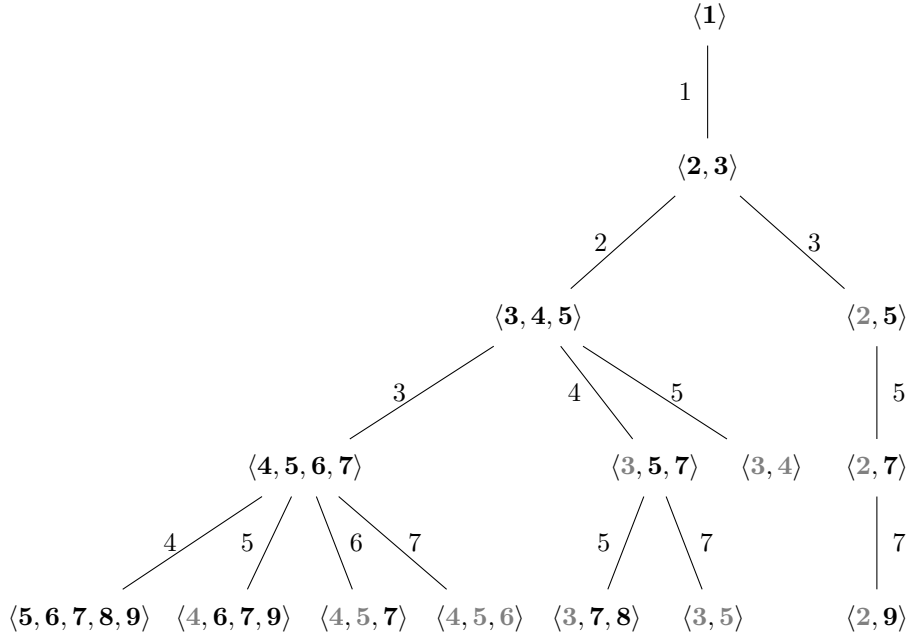
\begin{figure}[h!]
\[
\begin{tikzpicture}[level distance=2cm, inner sep=2mm,level/.style={sibling distance=1.55cm}]
    \node {$\sgnode{1}$}
    child {node {$\sgnode{2,3}$}
      child [sibling distance=4.5cm] {node {$\sgnode{3,4,5}$}
        child {node {$\sgnode{4,5,6,7}$}
          child [sibling distance=2cm]{node {$\sgnode{5,6,7,8,9}$}  edge from parent node [left] {4}}
          child [sibling distance=1.9cm]{node {$\sgnode{\gr 4,6,7,9}$}  edge from parent node [left] {5}}
          child {node {$\sgnode{\gr 4,\gr 5,7}$}  edge from parent node [right] {6}}
          child {node {$\sgnode{\gr 4,\gr 5,\gr 6}$}  edge from parent node [right] {7}}
          edge from parent node [left] {3}
        }
        child{{} edge from parent[draw=none]}
        child{{} edge from parent[draw=none]}
        child {node {$\sgnode{\gr 3,5,7}$}
          child {node {$\sgnode{\gr 3,7,8}$} edge from parent node [left] {5}}
          child {node {$\sgnode{\gr 3,\gr 5}$} edge from parent node [right] {7}}
          edge from parent node [left] {4}
        }
        child {node {$\sgnode{\gr 3,\gr 4}$}
          edge from parent node [right] {5}
        }
        edge from parent node [left] {2}
      }
      child [sibling distance=4.5cm] {node {$\sgnode{\gr2,5}$}
        child {node {$\sgnode{\gr2,7}$}
          child {node {$\sgnode{\gr2,9}$}
            edge from parent node [right] {7}
          }
          edge from parent node [right] {5}
        }
        edge from parent node [right] {3}
      }
      edge from parent node [left] {1}
    }
    ;
  \end{tikzpicture}
  \]
\caption{The part of the semigroup tree consisting of semigroups of genus at most $4$. A generator of a semigroup is in gray if it is not greater than $F(S)$. An edge between a semigroup $S$ and its child $S'$ is labeled by  $x$ if $S'$ is obtained from $S$ by removing $x$. This figure is taken from \cite{FromentinHivert}.}
\label{FigTree}
\end{figure}

In simple cases we can describe the set of children of a numerical semigroup explicitly.  For any positive integer $g$,
\[
\langle g+1,g+2,\ldots, 2g+1\rangle = \{0,g+1,g+2,\ldots\}
\]
is the \emph{ordinary semigroup} of genus $g$, which has Frobenius number $g$ and $g+1$ effective generators.  In \cite{BrasAmorosBounds}, Bras-Amor\'os computes the number of effective generators of each child of this semigroup and shows that the semigroup tree contains a subtree with exactly $2F_g$ nodes at level $g$.  Elizalde describes a more complicated subtree of the semigroup tree that leads to an improved lower bound \cite{Elizalde}.

\subsection{Strong generators and the proof of Theorem \ref{ZhaiTheorem}.}\ \\

One of the key ideas of Zhai's proof of Theorem \ref{ZhaiTheorem} is to divide the set of effective generators into two types.  Bulygin and Bras-Amor\'os define strong and weak effective generators in \cite{BrasAmorosBulygin}, and analyze the distribution of effective generators of each type in several classes of semigroups. For an effective generator $\lambda\in S$, the only element of $S$ that is not a minimal generator but may be a minimal generator of $S \setminus \{\lambda\}$ is $m(S) + \lambda$ \cite[Lemma 3]{BrasAmorosBounds}.

\begin{definition}
Let $S$ be a numerical semigroup.  An effective generator $\lambda$ is \emph{strong} if $m(S) + \lambda$ is a minimal generator of $S \setminus \{\lambda\}$.  An effective generator that is not strong is a \emph{weak generator}.

We say that $S'$ \emph{descends} from $S$ if $S' = S \setminus \{\lambda\}$ for some effective generator $\lambda$ of $S$.  This descent is \emph{strong} if $\lambda$ is a strong generator of $S$ and we say that $S'$ is a \emph{strongly descended} numerical semigroup.  A descent that is not strong is a \emph{weak descent}.  We say that $S''$ is a \emph{weak descendant} of $S$ if $S''$ is obtained from $S$ by a series of weak descents.
\end{definition}
%The following observation of Bras-Amor\'os motivates the definition of a strong generator.  Let $\lambda$ be a strong generator of $S$.  The only element of $S$ that is not a minimal generator but may be be a minimal generator of $S \setminus \lambda$, is $m(S) + \lambda$ \cite{BrasAmorosBounds}.

Every numerical semigroup is a weak descendant of a unique strongly descended numerical semigroup, where a strongly descended numerical semigroup is considered a weak descendant of itself.  Let $N_g(S)$ be the number of weak descendants of $S$ of genus $g$ and let $\mathcal{S}$ denote the set of strongly descended numerical semigroups.  Zhai studies $N(g)$ by analyzing the sum 
\begin{equation}\label{Eq1}
N(g) = \sum_{S\in \mathcal{S}} N_g(S).
\end{equation}
Since every weak descendant of $S$ comes from removing some subset of the $h(S)$ effective generators of $S$ we see that $N_g(S) \le \binom{h(S)}{g-g(S)}$.

Zhai further divides the sum in (\ref{Eq1}) into two pieces.  If $h(S) + g(S) < g$, then the bound from the previous paragraph shows that $N_g(S) = 0$.  Let $\mathcal{S}_2$ denote the set of strongly descended semigroups $S$ such that $h(S) + g(S) \ge g$ and $g(S) - h(S) < g/3$, and let $\mathcal{S}_3$ denote the set of strongly descended semigroups $S$ such that $h(S) + g(S) \ge g$ and $g(S) - h(S) \ge g/3$. We see that 
\begin{equation}\label{Eq2}
N(g) = \sum_{S\in \mathcal{S}_2} N_g(S) + \sum_{S\in \mathcal{S}_3} N_g(S).
\end{equation}
We write $N_2(g)$ for the first sum in (\ref{Eq2}) and $N_3(g)$ for the second.  The following proposition connects $N_2(g)$ to Zhao's results from Section \ref{Fibonacci}.
\begin{prop}[Zhai]
\begin{enumerate}
\item Every semigroup $S$ in $\mathcal{S}_2$ satisfies $F(S) < 2m(S)$.
\item If $S'$ is a weak descendant of a numerical semigroup in $\mathcal{S}_2$ then $F(S') < 3m(S')$. 
\end{enumerate}
\end{prop}

The second statement immediately implies that $N_2(g) \le t(g)$.  The heart of the proof of Theorem \ref{ZhaiTheorem} is to show that $N_2(g) = O(\phi^g)$ and that $N_3(g) = o(\phi^g)$, proving Conjecture \ref{ZhaoConjecture}.  Both of these bounds rely on the following estimate related to the set of strongly descended numerical semigroups with given multiplicity and Frobenius number. The proof of this result involves an intricate analysis of the semigroup tree.
\begin{lemma}[Zhai]\label{ZhaiLemma}
Let $\mathcal{S}(m,F)$ be the set of strongly descended numerical semigroups with multiplicity $m$ and Frobenius number $F$.  Then,
\[
\sum_{S \in \mathcal{S}(m,F)} \phi^{-(g(S)-h(S))} \le 5(F-m+2)\left(\frac{1.618}{\phi} \right)^{F-m-1}.
\]
\end{lemma}

\section{Properties of a ``typical'' numerical semigroup.}

Theorem \ref{ZhaiTheorem} implies that as $g$ goes to infinity, the average number of children of a numerical semigroup of genus $g$, or equivalently, the average number of effective generators, approaches $\phi$.  Let $t(g,h)$ be the number of numerical semigroups of genus $g$ with $h$ effective generators.  We would like to understand how $t(g,h)$ increases when $h$ is fixed and $g$ grows.  In particular, since $N(g+1) = \sum_h t(g,h)  h$, the only way for $N(g) > N(g+1)$ is if there are ``too many'' semigroups of genus $g$ that have no effective generators.

Lynnelle Ye, then an undergraduate, made the following conjecture in \cite{Ye}.
\begin{conjecture}[Ye]
For all $h \ge 0$, 
\[
\lim_{g \to \infty} \frac{t(g,h)}{N(g)} = \frac{1}{\phi^{h+2}}.
\]
\end{conjecture}
A stronger form of the conjecture was proven by Evan O'Dorney, also an undergraduate at the time, where the sum is taken over all values of $h$ instead of just considering an individual fixed value \cite{ODorney}.
\begin{theorem}[O'Dorney]\label{ODorneyTheorem}
We have
\[
\sum_{h \ge 0} \left| t(g,h) - S \phi^{g-(h+2)}\right| = o(\phi^g),
\]
where $S$ is the constant from Theorem \ref{ZhaiTheorem}.
\end{theorem}

O'Dorney also gives a more direct interpretation of the constant $S$ than had appeared previously.  Let $s(g,h)$ denote the number of strongly descended numerical semigroups with genus $g$ and $h$ effective generators. For any positive integer $n$ let $r(n) = s(2n+1,n+1)$ and define $r(-1) = r(0) = 1$.
\begin{prop}[O'Dorney]
We have
\[
S = \frac{\phi^2}{\sqrt{5}} \sum_{k\ge -1} r(k) \phi^{1-k},
\]
where $S$ is the constant from Theorem \ref{ZhaiTheorem}.
\end{prop}

Theorem \ref{ODorneyTheorem} shows that asymptotically the proportion of numerical semigroups of genus $g$ that have no children is $\phi^{-2} \approx .382$.  As long as $t(g,0) - \phi^{-2} N(g)$ is not ``too large'' for any particular value of $g$ we should be able to show that $N(g+1) \ge N(g)$.  This leads to questions about the error terms in Theorems \ref{ZhaiTheorem} and \ref{ODorneyTheorem}.  Unfortunately, these error terms are not currently in a useful form for explicit computations, which is due to the fact that $\frac{1.618}{\phi}$, which appears in Lemma \ref{ZhaiLemma}, is just barely less than $1$.  This difficulty with error terms makes it unclear how to use these results to prove Conjecture \ref{IncreasingConjecture}.  Ye uses the tools described above to prove a weaker statement \cite{Ye}.
\begin{prop}[Ye]
Let $\mathcal{S}(g)$ be the number of strongly descended numerical semigroups with genus $g$ and $\mathcal{N}_g$ denote the set of semigroups with genus $g$.  Then for all $g \ge 0$, we have 
\[
N(g+2) = N(g+1) - N(g) + \mathcal{S}(g+1) + 1 +  \sum_{S \in \mathcal{N}_g} \binom{h(S)-1}{2}.
\]
As a consequence, $N(g+2) \ge N(g+1) - N(g)$ for all $g$.
\end{prop}

One interpretation of Conjecture \ref{ZhaoConjecture} is that the number of numerical semigroups $S$ of genus $g$ with $F(S) > 3m(S)$ is $o(\phi^g)$, that a generically chosen numerical semigroup has Frobenius number at most $3$ times its multiplicity.  Starting from Zhai's proof of this conjecture, using Zhao's characterization of numerical semigroups with $F(S) < 3m(S)$ from Section \ref{Fibonacci}, Kaplan and Ye show that almost all numerical semigroups have their Frobenius number in a much smaller range \cite{KaplanYe}.
\begin{prop}[Kaplan, Ye]
Let $\epsilon > 0$ and $A_\epsilon(g)$ be the number of numerical semigroups $S$ with genus $g$ and $(2-\epsilon)m(S) < F(S) < (2+\epsilon)m(S)$.  Then $\lim_{g\to \infty} \frac{A_{\epsilon}(g)}{N(g)} = 1$.
\end{prop}
\noindent They also give a similar counting result for the numerical semigroups with $m/g$ in a particular range.
\begin{prop}[Kaplan, Ye]
Let $\epsilon > 0$ and $\gamma = \frac{5+\sqrt{5}}{10}$.  Let $\Phi_{\epsilon}(g)$ be the number of numerical semigroups with genus $g$ and $(\gamma-\epsilon)g < m(S) < (\gamma + \epsilon) g$.  Then $\lim_{g\to\infty} \frac{\Phi_{\epsilon}(g)}{N(g)} = 1$.
\end{prop}
\noindent This result shows that almost all numerical semigroups satisfy $2g(S) < 3m(S)$, which relates back to the discussion of Theorem \ref{2g3mTheorem}.

\section{Further questions.}

\subsection{Computing the semigroup tree.}

Theorem \ref{ZhaiTheorem} shows that $N(g) > N(g+1)$ can only hold for finitely many $g$.  By giving explicit error terms for several of the estimates in \cite{Zhai}, it should be possible to find an upper bound for the largest $g$ for which this is possible.  This would reduce Conjecture \ref{IncreasingConjecture} to a finite computation.

Bras-Amor\'os computed $N(g)$ for $g\le 50$ and states in \cite{BrasAmorosFibonacci} that the computation for $g=50$ took $19$ days.  Fromentin and Hivert use a massively improved algorithm for computing the semigroup tree, utilizing depth first rather than breadth first search along with several specific technical optimizations, to compute $N(g)$ for $g\le 67$.  More recently, Bras-Amor\'os and Fern\'andez-Gonz\'alez have suggested a new algorithm based on seeds, which can be thought of as a generalization of the notion of strong and weak effective generators \cite{BrasAmorosFernandezGonzalez}.

Since computing the full tree of all semigroups of bounded genus seems so computationally difficult, it may appear that computing $N(g)$ for $g$ large is hopeless.  There is another approach to this problem using a bijection between numerical semigroups of fixed multiplicity and integer points in a certain rational polyhedral cone.  

\subsection{Numerical semigroups and integer points in polytopes.}
Let $S$ be a numerical semigroup with multiplicity $m$ and recall from Section \ref{Fibonacci} that the Kunz coordinate vector, or Ap\'ery tuple, of $S$ is $(k_1,\ldots, k_{m-1})$ where the $k_i$ are positive integers defined so that $\Ap(m,S) = \{0,k_1 m+1,\ldots, k_{m-1}m +m-1\}$.  In this way, every numerical semigroup of multiplicity $m$ corresponds uniquely to an integer point in $\R^{m-1}$ \cite{Kunz}. The following result shows that the tuples of $m-1$ positive integers that arise as the Ap\'ery tuple of a numerical semigroup of multiplicity $m$ are exactly the integer points of a rational polyhedral cone \cite{BGGR}.

\begin{prop}[Branco, Garc\'ia-Garc\'ia, Garc\'ia-S\'anchez, Rosales] \label{NMG Ineq}
Consider the following set of inequalities:
\begin{eqnarray*}
x_i \ge 1 &\ \ \  &\ \ \  \text{for all}\ \ \ \ i\in \{1,\ldots, m-1\}, \\
x_i + x_j \ge x_{i+j} &\ \ \  &\ \ \ \text{for all}\ \ \ \ 1 \le i \le j \le m-1,\ i+j \le m-1,\\
x_i +x_j + 1 \ge x_{i+j-m} &\ \ \  &\ \ \ \text{for all}\ \ \ \ \ 1 \le i \le j \le m-1,\ i+j > m.
\end{eqnarray*}

There is a one-to-one correspondence between solutions $(k_1,\ldots, k_{m-1})$ to the above inequalities, with each $k_i \in \Z$, and the Ap\'ery tuples of numerical semigroups with multiplicity $m$.  

If we add the condition that $\sum_{i = 1}^{m-1} k_i = g$, then there is a one-to-one correspondence between solutions $(k_1,\ldots, k_{m-1})$ to the above inequalities, with each $k_i \in \Z$, and the Ap\'ery tuples of numerical semigroups with multiplicity $m$ and genus $g$. 
\end{prop}

The problem of computing $N(m,g)$ is thus reduced to the problem of counting the integer points in an $(m-2)$-dimensional rational polytope.  Blanco, Garc\'ia-S\'anchez, and Puerto use this characterization to show that for fixed $m,\ N(m,g)$ can be computed in polynomial time \cite{BlancoGarciaSanchezPuerto}.  Since $N(g) = \sum_{m \le g+1} N(m,g)$ we see that $N(g)$ can also be computed in polynomial time.  So far, this result has been of more theoretical than practical interest.  For instance, it is unclear how one would compute $N(m,g)$ for $m = 33$ and $g = 50$ with these techniques.

Several patterns seem to emerge when looking closely at Figure \ref{Fig4}.  For example, it is not difficult to guess that $N(3,g) = \left\lceil \frac{g+1}{3}\right\rceil$.  The authors of \cite{BlancoGarciaSanchezPuerto} prove this and a similar, but more complicated, formula for $N(4,g)$.  Kaplan uses results from Ehrhart theory to show that for fixed $m,\ N(m,g)$ is eventually given by a quasipolynomial of degree $m-2$ \cite{Kaplan}.  He also conjectures the following.
\begin{conjecture}[Kaplan]
For any $m \ge 2,\ N(m,g) \le N(m,g+1)$.
\end{conjecture}
This conjecture clearly implies Conjecture \ref{IncreasingConjecture}.  One proof strategy involves a detailed understanding of the polytopes defined by Proposition \ref{NMG Ineq}.

\subsection{The ordinarization transform.}

Bras-Amor\'os suggests another approach to Conjecture \ref{IncreasingConjecture} via the \emph{ordinarization transform} \cite{BrasAmorosOrdinarization} .  Let $S$ be a numerical semigroup with genus $g$, multiplicity $m$, and Frobenius number $F$.  It is easy to check that $(S \cup \{F\}) \setminus \{m\}$ is also a numerical semigroup of genus $g$.  Repeating this process, removing the multiplicity and adding the Frobenius number, gives a new numerical semigroup of genus $g$ each time, until we reach the ordinary semigroup of genus $g$, which was defined in Section \ref{SemigroupTree}.  The number of steps needed to reach the ordinary semigroup is called the \emph{ordinarization number} of $S$.  The ordinarization transform organizes all genus $g$ semigroups into a tree $\mathcal{T}_g$ where the root is given by the ordinary semigroup.  For a picture of $\mathcal{T}_6$ see \cite[Figure 1]{BrasAmorosOrdinarization}.

Bras-Amor\'os proposes the following conjecture, which implies Conjecture \ref{IncreasingConjecture}, about the number of nodes at level $r$ in the tree $\mathcal{T}_g$ compared to the number of nodes at level $r$ in $\mathcal{T}_{g+1}$. 
\begin{conjecture}[Bras-Amor\'os]
Let $n_{g,r}$ denote the number of numerical semigroups with genus $g$ and ordinarization number $r$.  For each genus $g \in \N_0$ and each ordinarization number $r\in \N_0,\ n_{g,r} \le n_{g+1,r}$.
\end{conjecture}
\noindent See \cite{BrasAmorosOrdinarization} for computational evidence and proofs of some cases.

\subsection{Counting by other invariants.}

The focus of this paper is on counting numerical semigroups ordered by genus, but we can also ask what happens for a different choice of ordering.  For example, let $\operatorname{ns}(F)$ denote the number of numerical semigroups with Frobenius number $F$.  One can compute that $\operatorname{ns}(5) = 5 > \operatorname{ns}(6) = 4$ and $\operatorname{ns}(31) = 70854> \operatorname{ns}(32) = 68681$ \cite[Table on page 11]{GarciaSanchezRosales}, so the analogue of Conjecture \ref{IncreasingConjecture} does not hold.  Restricting to values of $F$ of the same parity helps to clarify the overall growth rate. Backelin proves a type of analogue of Theorem \ref{ZhaiTheorem} in \cite{Backelin}.

\begin{theorem}[Backelin]
We have that 
\[
\lim_{F \to \infty \atop F\ \text{odd}} 2^{-F/2} \operatorname{ns}(F)\ \ \ \ \  \text{ and }\ \ \ \ \   \lim_{F \to \infty \atop F\ \text{even}} 2^{-F/2}\operatorname{ns}(F)
\] 
exist and are nonzero.
\end{theorem}
\noindent The values of these limits are not known. It is unclear if we should expect them to be the same.

The \emph{weight} of a numerical semigroup $S$ of genus $g$, denoted $w(S)$, is equal to the sum of the gaps of $S$ minus $g(g+1)/2$, or equivalently, 
\[
w(S) = \sum_{i=1}^g (l_i - i),
\] 
where $l_1,\ldots, l_g$ are the gaps of $S$.  For example, the ordinary numerical semigroup of genus $g$ has weight $0$.  The weight plays an important role in the connection between numerical semigroups and algebraic curves, which we discuss in the next section.  Since there are infinitely many numerical semigroups of weight $0$ it does not make sense to try to count semigroups ordered by weight.  

\begin{figure}
\begin{center}
\begin{tikzpicture}[scale = 0.6]

				\begin{scope}[xshift = 16cm]
					\foreach \x in {0,1,2}
						{ \draw (\x,3) rectangle ++(1,1);}
					\foreach \x in {0,1,2}
						{ \draw (\x,2) rectangle ++(1,1);}
					\foreach \x in {0}
						{ \draw (\x,1) rectangle ++(1,1);}
					\foreach \x in {0}
						{ \draw (\x,0) rectangle ++(1,1);}
					
					\begin{scope}[blue, ultra thick]
						\filldraw (0,0) circle [radius = 0.1];
						\draw[font = \footnotesize]
							(0,0) --
							node [below] {0} ++(1,0) --
							node [right] {1} ++(0,1) -- 
							node [right] {2} ++(0,1) --
							node [below right] {3} ++(1,0) --
							node [below] {4} ++(1,0) --
							node [right] {5} ++(0,1) --
							node [right] {6} ++(0,1.05);
					\end{scope}
					
				\end{scope}

			\end{tikzpicture}
		\end{center}
		\caption{The partition corresponding to $\langle 3,4\rangle = \{0,3,4,7,8,\ldots\}$.}
		\label{FigPartition}
\end{figure}
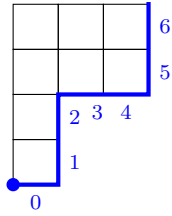

Bras-Amor\'os and de Mier show how the \emph{enumeration} of a numerical semigroup gives rise to a Dyck path inside of a square \cite{BrasAmorosdeMier}.  Given a numerical semigroup $S$ define $\tau(S)$ as the path with origin $(0,0)$ and steps $e(i)$ given by
\[
e(i) = \begin{cases}
\rightarrow & \text{ if } i \in S,\\
\uparrow & \text{ if } i \not\in S,
\end{cases}\ \ \ \ \ \ \ \ \text{ for } 0 \le i \le 2g(S).
\]
This is a slight variation of the construction in \cite{BrasAmorosdeMier}, since we include an initial step to the right for $e(0)$.  Taking this path, together with the $x$-axis and the line $y = g$ gives the Ferrers diagram of a partition of size equal to $w(S) + g(S)$.  Figure \ref{FigPartition} gives the partition corresponding to $\langle 3,4\rangle$.

This connection between numerical semigroups and partitions is the subject of a paper by Kaplan, Hannah Constantin, and Benjamin Houston-Edwards, both undergraduates at the time \cite{ConstantinHoustonEdwardsKaplan}.  Very little is known about the number of numerical semigroups with fixed weight plus genus.  For further results counting numerical semigroups in certain families by counting lattice paths with special properties, see \cite{HellusWaldi, KunzWaldi}.

\subsection{Numerical semigroups and algebraic curves.}

Let $C$ be a smooth projective complex curve.  At every point $p \in C$ we consider the set of orders of poles of rational functions that are regular away from $p$.  These pole orders form a numerical semigroup called the \emph{Weierstrass semigroup} at $p$.  The Weierstrass \emph{L\"uckensatz} (or Gap Theorem) says that the genus of the semigroup is equal to the genus of the curve $C$.  At almost all points of $C$ the associated semigroup is the semigroup $\{0,g+1,g+2,\ldots\}$, which is the motivation for calling this semigroup ordinary \cite{BrasAmorosIEEE}.  Any point for which this is not the case is called a \emph{Weierstrass point}.  The multiset of Weierstrass semigroups at these points reflects the geometry of the curve.  For instance, a genus $g$ curve is hyperelliptic if and only if it has a point with Weierstrass semigroup $\langle 2, 2g+1\rangle$.  See del Centina's article \cite{delCentina} for a nice historical overview of the subject.

In the late 19\textsuperscript{th} century, Hurwitz asked for a characterization of the numerical semigroups that occur as the Weierstrass semigroup of some point on some curve.  Buchweitz proved that not every semigroup $S$ arises, giving a criterion in terms of sumsets of the gaps of $S$ \cite{Buchweitz}.  Kaplan and Ye show that the number of genus $g$ semigroups failing this criterion is $o(\phi^g)$ \cite{KaplanYe}.  We do not know whether a positive proportion of semigroups occur as Weierstrass semigroups.

We would like to understand not only whether a semigroup occurs as a Weierstrass semigroup, but also the dimension of the space of genus $g$ curves with such a Weierstrass semigroup.  With this problem in mind, Pflueger defines the following variation of the weight of a semigroup \cite{Pflueger}.

\begin{definition}
Let $S$ be a numerical semigroup with gaps $l_1,\ldots, l_g$ and minimal generating set $n_1,\ldots, n_e$.  The \emph{effective weight} of $S$ is 
\[
\operatorname{ewt}(S) = \sum_{i=1}^g \#\{j\ |\ n_j < l_i\}.
\]
That is, the effective weight is the sum over all gaps, of the number of minimal generators less than that gap.
\end{definition}

Computing this quantity for all semigroups with genus $g \le 50$ leads to the following purely combinatorial problem \cite{Pflueger}.
\begin{conjecture}[Pflueger]
Let $S$ be a numerical semigroup of genus $g$.  Then 
\[
\operatorname{ewt}(S) \le \left\lfloor \frac{(g+1)^2}{8}\right\rfloor.
\]
\end{conjecture}

\subsection{The Wilf conjecture.}

In a 1978 article in this \textsc{Monthly}, Wilf proposed the following problem, which has become one of the most studied questions in the theory of numerical semigroups \cite{Wilf}.  It was originally phrased as a question, but has come to be known as Wilf's conjecture, as most authors seem to believe that it is true.

\begin{conjecture}[Wilf]\label{WilfConj}
Let $S$ be a numerical semigroup with embedding dimension $e$, Frobenius number $F$, and $|S \cap [0,F]| = n$.  Then 
\[
F+1  \le ne.
\]
\end{conjecture}
This conjecture has been verified for all semigroups with $g \le 60$ by Fromentin and Hivert \cite{FromentinHivert}.  Results on this problem due to Alex Zhai \cite{ZhaiWilf}, and to Alessio Sammartano and Alessio Moscariello \cite{MoscarielloSammartano, Sammartano}, were started as undergraduate research.  There are many special cases known, but Conjecture \ref{WilfConj} is open in general.  See for example, the work of Dobbs and Matthews \cite{DobbsMatthews}, Kaplan \cite{Kaplan}, and the recent papers of Eliahou \cite{Eliahou}, and Delgado \cite{Delgado}.

\subsection{Higher dimensions.}

The main focus of this article is counting submonoids of $\N_0$ ordered by the size of their complement.  There is a natural higher-dimensional version of this problem.  For fixed $d$, let $N_d(g)$ denote the number of submonoids $S$ of $\N_0^d$ for which $|\N_0^d\setminus S| = g$.  Failla, Peterson, and Utano study $N_d(g)$ in \cite{FaillaPetersonUtano}, giving lower bounds coming from some special classes of submonoids and investigating analogues of the semigroup tree.  A higher-dimensional analogue of Wilf's conjecture is proposed in \cite{GarciaGarciaMarinAragonVigneronTenorio}. 

For fixed $d$, how does $N_d(g)$ grow?  This problem is completely open.  Even for $d = 2$ we do not have a conjecture for the growth rate.

\section{Acknowledgment.}
The author thanks Scott Chapman for introducing him to numerical semigroups.  He thanks Joe Gallian for encouraging him to write this article and for comments on a draft of this paper.  He thanks Nathan Pflueger for providing data and for helpful conversations.  He also thanks the referees for valuable suggestions.

\end{document}